\documentclass[12pt,a4paper]{article}
\usepackage{indentfirst}
\usepackage{amsmath,amssymb,amsfonts,amsthm,graphics}
\usepackage{makeidx}
\newtheorem{thm}{Theorem}[section]
\newtheorem{definition} {Definition}[section]
\newtheorem{pro}{Proposition}[section]
\newtheorem{lem}{Lemma}[section]

\newtheorem{conj}{Conjecture}[section]

\def\pf{\noindent {\it Proof.} }

\begin{document}
\title{\Large\bf Sharp upper bound for the rainbow connection\\ numbers of
$2$-connected graphs\footnote{Supported by NSFC No.11071130.}}
\author{\small Xueliang Li, Sujuan Liu\\
\small Center for Combinatorics and LPMC-TJKLC\\
\small Nankai University, Tianjin 300071, China\\
\small lxl@nankai.edu.cn;
sjliu0529@126.com}
\date{}
\maketitle
\begin{abstract}
An edge-colored graph $G$, where adjacent edges may be colored the
same, is rainbow connected if any two vertices of $G$ are connected
by a path whose edges have distinct colors. The rainbow connection
number $rc(G)$ of a connected graph $G$ is the smallest number of
colors that are needed in order to make $G$ rainbow connected. In
this paper, we give a sharp upper bound that
$rc(G)\leq\lceil\frac{n}{2}\rceil$ for any $2$-connected graph $G$
of order $n$, which improves the results of Caro et al. to best
possible.

\end{abstract}

\section{Introduction}

All graphs considered in this paper are simple, finite and
undirected. An $edge$-$coloring$ of a graph $G$ is a function from
the edge set of $G$ to the set of natural numbers. A path in an
edge-colored graph $G$ is a $rainbow\ path$ if no two edges of this
path are colored the same. An edge-colored graph is $rainbow\
connected$ if every pair of vertices is connected by at least one
rainbow path. The $rainbow\ connection\ number$ of a connected graph
$G$, denoted by $rc(G)$, is the smallest number of colors that are
needed to rainbow color the graph $G$. We call a rainbow coloring of
$G$ with $k$ colors a {\it $k$-rainbow coloring}.

The concept of a rainbow coloring was introduced by Chartrand et al.
in \cite{char}. The rainbow connection numbers of several graph
classes have been obtained. It is well known that a cycle with $n$
vertices has a rainbow connection number $\lceil\frac{n}{2}\rceil$,
$rc(G)=n-1$ if and only if $G$ is a tree of order $n \ (\geq2)$, and
$rc(G)=1$ if and only if $G$ is a complete graph of order $n \ (\geq
2)$. In \cite{chak}, Chakraborty et al. gave the following result
about the rainbow connection number.

\begin{thm} \cite{chak} Given a graph $G$,
deciding if $rc(G)=2$ is NP-Complete. In particular, computing
$rc(G)$ is NP-Hard. \end{thm}

However, many upper bounds of the rainbow connection number have
been given. For a $2$-connected graph, Caro et al. proved the
following two results.

\begin{pro} \cite{Caro} If $G$ is a
$2$-connected graph with $n$ vertices, then $rc(G)\leq\frac{2n}{3}$.
\end{pro}

\begin{thm} \cite{Caro} If $G$ is a $2$-connected graph on $n$ vertices,
then $rc(G)\leq\frac{n}{2}+O(n^{\frac{1}{2}})$.
\end{thm}

One can see that both the above bounds are much greater than
$\lceil\frac{n}{2}\rceil$. However, experience tells us that the
best bound should be $\lceil\frac{n}{2}\rceil$. This paper is to
give it a proof. Before proceeding, we need the following notation
and terminology.

A $separation$ of a connected graph is a decomposition of the graph
into two nonempty connected subgraphs which have just one vertex in
common. This common vertex is called a $separating$ $vertex$ of the
graph. A graph is $nonseparable$ if it is connected and has no
separating vertices; otherwise, it is separable. If a graph $G$ has
at least $3$ vertices but no loops, then $G$ is nonseparable if and
only if $G$ is $2$-connected.

Let $F$ be a subgraph of a graph $G$. An $ear$ of $F$ in $G$ is a
nontrivial path in $G$ whose ends are in $F$ but whose internal
vertices are not. A nested sequence of graphs is a sequence ($G_0,
G_1, \cdots, G_k$) of graphs such that $G_i\subset G_{i+1}, 0\leq i<
k$. An $ear$ $decomposition$ of a nonseparable graph $G$ is a nested
sequence ($G_0, G_1, \cdots, G_k$) of nonseparable subgraphs of $G$
such that: $(1)$ $G_0$ is a cycle; $(2)$ $G_i=G_{i-1}\bigcup P_i$,
where $P_i$ is an ear of $G_{i-1}$ in $G, 1\leq i\leq k$; $(3)$
$G_k=G$. It is clear that every $2$-connected graph has an ear
decomposition.

At the IWOCA workshop \cite{sch}, Hajo Broersma posed a question:
what happens with the value $rc(G)$ for graphs with higher
connectivity? Motivated by this question, we study the rainbow
connection number of a $2$-connected graph and give a sharp upper
bound that $rc(G)\leq\lceil\frac{n}{2}\rceil$ for any $2$-connected
graph of order $n$, which improves the results in \cite{Caro} to
best possible.

\section{Main results}

For convenience, we first introduce some new definitions.

\begin{definition} Let $c$ be a $k$-rainbow coloring
of a connected graph $G$. If a rainbow path $P$ in $G$ has length
$k$, we call $P$ a complete rainbow path; otherwise, it is a
noncomplete rainbow path. A rainbow coloring $c$ of $G$ is
noncomplete if for any vertex $u\in V(G)$, $G$ has {\bf at most one}
vertex $v$ such that all the rainbow paths between $u$ and $v$ are
complete; otherwise, it is complete.
\end{definition}

For a connected graph $G$, if a spanning subgraph has a
(noncomplete) $k$-rainbow coloring, then $G$ has a (noncomplete)
$k$-rainbow coloring. This simple fact will be used in the following
proofs.

\begin{lem} \label{lem1} Let $G$ be a Hamiltonian graph of order $n
\ (n\geq 3)$. Then $G$ has a noncomplete
$\lceil\frac{n}{2}\rceil$-rainbow coloring, i.e.,
$rc(G)\leq\lceil\frac{n}{2}\rceil$.
\end{lem}

\pf Since $G$ is a Hamiltonian graph, there is a Hamiltonian cycle
$C_n=v_1, v_2, \cdots , v_n, v_{n+1}$ $(=v_1)$ in $G$. Define the
edge-coloring $c$ of $C_n$ by $c(v_iv_{i+1})=i$ for $1\leq i\leq
\lceil\frac{n}{2}\rceil$ and
$c(v_iv_{i+1})=i-\lceil\frac{n}{2}\rceil$ if
$\lceil\frac{n}{2}\rceil+1\leq i\leq n$. It is clear that $c$ is a
$\lceil\frac{n}{2}\rceil$-rainbow coloring of $C_n$, and the
shortest path connecting any two vertices in $V(G)$ on $C_n$ is a
rainbow path. For any vertex $v_i \ (1\leq i\leq n)$, only the
antipodal vertex of $v_i$ has no noncomplete rainbow path to $v_i$
if $n$ is even. Every pair of vertices in $G$ has a noncomplete
rainbow path if $n$ is odd. Hence the rainbow coloring $c$ of $C_n$
is noncomplete. Since $C_n$ is a spanning subgraph of $G$, $G$ has a
noncomplete $\lceil\frac{n}{2}\rceil$-rainbow coloring.
$\blacksquare$

Let $G$ be a 2-connected non-Hamiltonian graph of order $n \ (n\geq
4)$. Then $G$ must have an even cycle. In fact, since $G$ is
2-connected, $G$ must have a cycle $C$. If $C$ is an even cycle, we
are done. Otherwise, $C$ is a odd cycle, we then choose an ear $P$
of $C$ such that $V(C)\bigcap V(P)=\{a, b\}$. Since the lengths of
the two segments between $a, b$ on $C$ have different parities, $P$
joining with one of the two segments forms an even cycle. Then,
starting from an even cycle $G_0$, there exists a nonincreasing ear
decomposition $(G_0, G_1,$ $ \cdots, G_t, G_{t+1},$ $ \cdots, G_k)$
of $G$, such that $G_i=G_{i-1}\bigcup P_i \ (1\leq i\leq k)$ and
$P_i$ is a longest ear of $G_{i-1}$, i.e., $\ell(P_1)\geq \ell(P_2)
\geq \cdots \geq \ell(P_k)$. Suppose that $V(P_i) \bigcap
V(G_{i-1})=\{a_i, b_i\} \ (1\leq i\leq k)$. We call the distinct
vertices $a_i, b_i \ (1\leq i\leq k)$ the $foot$ $vertices$ of the
ear $P_i$. Without loss of generality, suppose that $\ell(P_t)\geq
2$ and $\ell(P_{t+1})= \ \cdots \ = \ell(P_k)=1$. So $G_t$ is a
2-connected spanning subgraph of $G$. Since $G$ is a non-Hamiltonian
graph, we have $t\geq 1$. Denote the order of $G_i \ (0\leq i\leq
k)$ by $n_i$. All these notations will be used in the sequel.

\begin{lem} \label{lem2} Let $G$ be a 2-connected non-Hamiltonian
graph of order $n \ (n\geq 4)$. If $G$ has at most one ear with
length 2 in the nonincreasing ear decomposition, then $G$ has a
noncomplete $\lceil\frac{n}{2}\rceil$-rainbow coloring, i.e.,
$rc(G)\leq\lceil\frac{n}{2}\rceil$.
\end{lem}

\pf Since $G_t \ (t\geq1)$ in the nonincreasing ear decomposition is
a 2-connected spanning subgraph of $G$, it only needs to show that
$G_t$ has a noncomplete $\lceil\frac{n_t}{2}\rceil$-rainbow coloring.
We will apply induction on $t$.

First, consider the case that $t=1$. Let $G$ be a 2-connected
non-Hamiltonian graph with $t=1$ in the nonincreasing ear
decomposition. The spanning subgraph $G_1=G_0\bigcup P_1$ of $G$
consists of an even cycle $G_0$ and an ear $P_1$ of $G_0$. Without
loss of generality, suppose that $G_0=v_1, v_2, \cdots , v_{2k},
v_{2k+1} \ (=v_1)$ where $k\geq2$. We color the edges of $G_0$ with
$k$ colors. Define the edge-coloring $c_0$ of $G_0$ by $c_0
(v_iv_{i+1})=i$ for $1\leq i\leq k$ and $c_0(v_iv_{i+1})=i-k$ if
$k+1\leq i\leq 2k$. From the proof of Lemma \ref{lem1}, the coloring
$c_0$ is a noncomplete $k$-rainbow coloring of $G_0$. Now consider
$G_1$ according to the parity of $\ell(P_1)$. If $\ell(P_1)$ is
even, then $n_1$ is odd and color the edges of $P_1$ with
$\frac{\ell(P_1)}{2}$ new colors. In the first $\frac{\ell(P_1)}{2}$
edges of $P_1$ the colors are all distinct, and the same ordering of
colors is repeated in the last $\frac{\ell(P_1)}{2}$ edges of $P_1$.
It is easy to verify that the obtained coloring $c_1$ of $G_1$ is a
noncomplete $\lceil\frac{n_1}{2}\rceil$-rainbow coloring and that
for any pair of vertices in $G$, there exists a noncomplete rainbow
path connecting them. If $\ell(P_1)$ is odd, then $n_1$ is even and
color the edges of $P_1$ with $\frac{\ell(P_1)-1}{2}$ new colors.
The middle edge of $P_1$ receives any color that already appeared in
$G_0$. The first $\frac{\ell(P_1)-1}{2}$ edges of $P_1$ all receive
distinct new colors and in the last $\frac{\ell(P_1)-1}{2}$ edges of
$P_1$ this coloring is repeated in the same order. It is easy to
verify that the obtained coloring $c_1$ of $G_1$ is a noncomplete
$\lceil\frac{n_1}{2}\rceil$-rainbow coloring.

Let $G$ be a 2-connected non-Hamiltonian graph with $t\geq 2$ in the
nonincreasing ear decomposition. Assume that the subgraph $G_i \
(1\leq i \leq t-1)$ has a noncomplete
$\lceil\frac{n_{i}}{2}\rceil$-rainbow coloring $c_i$ and when
$n_{i}$ is odd, any pair of vertices have a noncomplete rainbow
path. We distinguish the following three cases.

\noindent {\bf Case 1.}  $\ell(P_t) \ (\geq 3)$ is odd.

Suppose that $P_t=v_0(=a_t), v_1, \cdots , v_r, v_{r+1}, \cdots ,
v_{2r}, v_{2r+1}(=b_t)$ where $r\geq1$. We color the edges of $P_t$
with $r$ new colors to obtain a noncomplete coloring $c_t$ of $G_t$.
Define the edge-coloring of $P_t$ by $c(v_{i-1}v_i)=x_i \ (1\leq i
\leq r)$, $c(v_rv_{r+1})=x$ and $c(v_{i-1}v_i)=x_{i-r-1} \ (r+2\leq
i \leq 2r+1)$, where $x_1, x_2, \cdots , x_r$ are new colors and $x$
is a color appeared in $G_{t-1}$. It is easy to check that the
obtained coloring $c_t$ of $G_t$ is a
$\lceil\frac{n_t}{2}\rceil$-rainbow coloring.

Now we show that $c_t$ is noncomplete. For any pair of vertices $u,
v\in V(G_{t-1})\times V(G_{t-1})$, the rainbow path $P$ from $u$ to
$v$ in $G_{t-1}$ is noncomplete in $G_t$, because the new colors
$x_1, x_2, \cdots , x_r \ (r\geq1)$ do not appear in $P$. For any
pair of vertices $u, v\in V(P_t)\times V(P_t)$, if there exists a
rainbow path $P$ from $u$ to $v$ on $P_t$, then $P$ is noncomplete
in $G_t$, since some color in $G_{t-1}$ does not appear in $P$; if
not, there exists a noncomplete rainbow path $P$ from $u$ to $v$
through some vertices of $G_{t-1}$ such that at least one new color
does not appear in $P$. For any pair of vertices $um v \in
V(G_{t-1}) \times (V(P_t)\backslash \{v_r,v_{r+1}\})$, there exists
a noncomplete rainbow path from $u$ to $v$ in which at least one new
color does not appear. If there exists a vertex all of whose rainbow
paths to $a_t$ (resp. $b_t$) in $G_{t-1}$ are complete, we denote
the vertex by $a_t^\prime$ (resp. $b_t^\prime$). For vertex $v_r$
(resp. $v_{r+1}$), only the vertex $a_t^\prime$ (resp. $b_t^\prime$)
possibly has no noncomplete rainbow path to $v_r$ (resp. $v_{r+1}$)
in $G_t$. So there possibly exist two pairs of vertices $a_t^\prime,
v_r$ and $b_t^\prime, v_r+1$ which have no noncomplete rainbow path.
Since $a_t^\prime, \ b_t^\prime$ are distinct in $G_{t-1}$, the
rainbow coloring $c_t$ is noncomplete. If $n_t$ is odd, then
$n_{t-1}$ is odd. By induction, $a_t^\prime, \ b_t^\prime$ do not
exist when $n_{t-1}$ is odd. Hence every pair of vertices have a
noncomplete rainbow path.

\noindent {\bf Case 2.} $\ell(P_t) \ (\geq 2)$ is even and $n_{t-1}$
is even.

In this case, $n_t$ is odd. Suppose that $P_t=v_0(=a_t), v_1, \cdots, v_r, v_{r+1},
\cdots, v_{2r-1}, v_{2r}(=b_t)$ where $r\geq1$. Define the
edge-coloring of $P_t$ by $c(v_{i-1}v_i)=x_i$ for $1\leq i\leq r$
and $c(v_{i-1}v_i)=x_{i-r}$ for $r+1\leq i\leq 2r$. It is
clear that the obtained coloring $c_t$ of $G_t$ is a
$\lceil\frac{n_t}{2}\rceil$-rainbow coloring.

Now we prove that $c_t$ is noncomplete. For any pair of vertices in
$V(G_{t-1})\times V(G_{t-1})$ or $V(P_t)\times V(P_t)$, there is a
noncomplete rainbow path connecting them in $G_t$, similar to the
Case $1$. For any pair of vertices $u\in V(G_{t-1}), v\in V(P_t) \
(v\neq v_r)$, there is a noncomplete rainbow path $P$ from $u$ to
$v$ such that at least one new color does not appear in $P$. For any
vertex $u\in V(G_{t-1})$, since the coloring $c_{t-1}$ is
noncomplete, $u$ has a noncomplete rainbow path $P^\prime$ in
$G_{t-1}$ to one of $a_t, \ b_t$ (say $a_t$). Then $P^\prime$
joining with $a_tP_tv_r$ is a noncomplete rainbow path from $u$ to
$v_r$ in $G_t$. Therefore, the rainbow coloring $c_t$ of $G_t$ is
noncomplete such that any pair of vertices has a noncomplete rainbow
path.

\noindent {\bf Case 3.} $\ell(P_t) \ (\geq 2)$ is even and $n_{t-1}$
is odd.

In this case, $n_t$ is even. We consider the following three subcases.

\noindent {\bf Subcase 3.1} $[V(P_t)\bigcap V(P_{t-1})]\backslash V(G_{t-2})=\emptyset$.

If $\ell(P_{t-1})$ is odd, let $G_{t-1}^\prime =G_{t-2}\bigcup P_t$
and $G_t=G_{t-1}^\prime \bigcup P_{t-1}$. By induction,
$G_{t-1}^\prime$ has a noncomplete
$\lceil\frac{n_{t-1}^\prime}{2}\rceil$-rainbow coloring
($n_{t-1}^\prime$ is the order of $G_{t-1}^\prime$). Similar to Case
1, we can obtain a noncomplete $\lceil\frac{n_t}{2}\rceil$-rainbow
coloring of $G_t$ from $G_{t-1}^\prime$ .

Suppose that $\ell(P_{t-1})$ is even. By induction, $G_{t-2}$ has a
noncomplete $\lceil\frac{n_{t-2}}{2}\rceil$-rainbow coloring
$c_{t-2}$. Suppose that $P_{t-1}=v_0(=a_{t-1}), v_1, \cdots , v_r,
v_{r+1},$ $\cdots, v_{2r-1}, v_{2r}(=b_{t-1})$ and $P_t=v_0^\prime
(=a_t), v_1^\prime, \cdots , v_s^\prime, v_{s+1}^\prime, \cdots ,
v_{2s-1}^\prime, v_{2s}^\prime(=b_t)$, where $r\geq 2, s\geq 1$.
Since $c_{t-2}$ is noncomplete and $a_i, b_i \ (1\leq i\leq k)$ are
two distinct vertices, then $a_{t-1}$ has a noncomplete rainbow path
$P^\prime$ to one of $a_t, b_t$ (say $a_t$) and $b_{t-1}$ has a
noncomplete rainbow path $P^{\prime\prime}$ to the other vertex.
Suppose that $x$ is the color in $G_{t-2}$ that does not appear in
$P^\prime$. Now color the edges of $P_{t-1}, P_t$ with $r+s-1$ new
colors and the color $x$. Define an edge-coloring of $P_{t-1}$ by
$c(v_{i-1}v_i)=x_i \ (1\leq i \leq r)$ and $c(v_{i-1}v_i)=x_{i-r} \
(r+1\leq i \leq 2r)$, where $x_1, x_2, \cdots , x_r$ are new colors.
And define an edge-coloring of $P_t$ by $c(v_{i-1}^\prime
v_i^\prime)=y_i \ (1\leq i \leq s-1), c(v_{s-1}^\prime v_s^\prime
)=x, c(v_s^\prime v_{s+1}^\prime )=x_1$ and $c(v_{i-1}^\prime
v_i^\prime )=y_{i-s-1} \ (s+2\leq i \leq 2s)$, where $y_1, y_2,
\cdots , y_{s-1}$ are new colors.

Similar to Case $2$, the obtained coloring $c_{t-1}$ of $G_{t-1}$ is
a noncomplete $\lceil\frac{n_{t-1}}{2}\rceil$-rainbow coloring such
that every pair of vertices have a noncomplete rainbow path. It is
obvious that $G_t$ is rainbow connected. The path $(v_s^\prime
P_ta_t)P^\prime(a_{t-1}P_{t-1}v_r)$ is a rainbow path from
$v_s^\prime$ to $v_r$ which is possibly complete. For any other pair
of vertices in $G_t$, there is a noncomplete rainbow path connecting
them. Hence the rainbow coloring $c_t$ of $G_t$ is noncomplete.

\noindent {\bf Subcase 3.2} $[V(P_t)\bigcap V(P_{t-1})]\backslash V(G_{t-2})=\{b_t\}$.

If $\ell(P_{t-1})$ is odd, suppose that $P_{t-1}=v_0(=a_{t-1}), v_1,
\cdots , v_r, v_{r+1}, \cdots , v_{2r}, v_{2r+1}(=b_{t-1})$. Since
$P_{t-1}$ is a longest ear of $G_{t-2}$ and $b_t\in
V(P_{t-1})\setminus V(G_{t-2})$, we have $r\geq 2$. Define an
edge-coloring of $P_{t-1}$ by $c(v_{i-1}v_i)=x_i \ (1\leq i \leq
r)$, $c(v_rv_{r+1})=x$ and $c(v_{i-1}v_i)=x_{i-r-1} \ (r+2\leq i
\leq 2r+1)$, where $x_1, x_2, \cdots , x_r$ are new colors and $x$
is a color appeared in $G_{t-2}$. Similar to Case $1$, the obtained
coloring $c_{t-1}$ of $G_{t-1}$ is a noncomplete
$\lceil\frac{n_{t-1}}{2}\rceil$-rainbow coloring such that every
pair of vertices have a noncomplete rainbow path. If $\ell(P_{t-1})$
is even, suppose that $P_{t-1}=v_0(=a_{t-1}), v_1, \cdots, v_r,
v_{r+1}, \cdots , v_{2r-1}, v_{2r}(=b_{t-1})$, where $r\geq2$.
Define an edge-coloring of $P_{t-1}$ by $c(v_{i-1}v_i)=x_i \ (1\leq
i \leq r)$, and $c(v_{i-1}v_i)=x_{i-r} \ (r+1\leq i \leq 2r)$, where
$x_1, x_2, \cdots , x_r$ are new colors. Similar to Case 2, the
obtained coloring $c_{t-1}$ of $G_{t-1}$ is a noncomplete
$\lceil\frac{n_{t-1}}{2}\rceil$-rainbow coloring such that every
pair of vertices have a noncomplete rainbow path.

Without loss of generality, assume that $b_t$ belongs to the first
half of $P_{t-1}$ and that $P_t=v_0^\prime (=a_t), v_1^\prime,
\cdots ,$ $ v_s^\prime, v_{s+1}^\prime, \cdots , v_{2s-1}^\prime,
v_{2s}^\prime(=b_t)$, where $s\geq1$. We color the edges of $P_t$
with $s-1$ new colors. Define an edge-coloring of $P_t$ by
$c(v_{i-1}^\prime v_i^\prime)= y_i \ (1\leq i \leq s-1)$,
$c(v_{s-1}^\prime v_s^\prime )=x_1, c(v_s^\prime v_{s+1}^\prime )
=y$ and $c(v_{i-1}^\prime v_i^\prime )=y_{i-s-1} \ (s+2\leq i \leq
2s)$, where $y_1, y_2, \cdots , y_{s-1}$ are new colors and the
color $y$ is different from color $x$ in $G_{t-2}$. It is easy to
verify that the obtained coloring $c_t$ of $G_t$ is a
$\lceil\frac{n_t}{2}\rceil$-rainbow coloring.

For any pair of vertices $v^\prime\in V(P_t)(v^\prime\neq v_s^\prime)$
and $v\in V(G_{t-1})$, there exists a noncomplete rainbow path $P$
connecting them since the path from $v^\prime$ to one foot
vertex of $P_t$ colored by new colors joining with the noncomplete
rainbow path from the foot vertex to $v$ in $V(G_{t-1})$ is a noncomplete
rainbow path from $v^\prime$ to $v$ in $G_t$. For $v_s^\prime$,
there is a noncomplete rainbow path from $v_s^\prime$ to any
vertex in $V(G_{t-2})\bigcup V(b_{t-1}P_{t-1}v_{r+2})$ through
edge $e=v_{s-1}^\prime v_s^\prime$; and a noncomplete rainbow
path from $v_s^\prime$ to any vertex in $V(a_{t-1}P_{t-1}v_{r+1})$
through $e=v_s^\prime v_{s+1}^\prime$. For any pair of vertices
in $V(P_t)\times V(P_t)$, there is a noncomplete rainbow path
connecting them obviously. Hence the rainbow coloring $c_t$ is noncomplete.

\noindent {\bf Subcase 3.3} $[V(P_t)\bigcap V(P_{t-1})]\backslash V(G_{t-2})=\{a_t, b_t\}$.

We can prove this subcase in a way similar to Subcase 3.2. Without
loss of generality, we can assume that $a_t=v_p(1\leq p\leq r-1)$
and $b_t=v_q(q\geq p+2)$. Color all the edges of $P_{t-1}$ and $P_t$
as in Subcase 3.2 but only the edge $e=v_{s-1}^\prime v_s^\prime$
which is colored by $x_{j+1}$ instead. The obtained coloring $c_t$
of $G_t$ is a noncomplete $\lceil\frac{n_t}{2}\rceil$-rainbow
coloring. $\blacksquare$

\begin{lem} \label{lem3} Let $G$ be a 2-connected non-Hamiltonian
graph of order $n \ (n\geq 4)$. If $G$ has at least $2$ ears of
length 2 in the nonincreaing ear decomposition, then
$rc(G)\leq\lceil\frac{n}{2}\rceil$.
\end{lem}

\pf We only need to prove that there exists a rainbow coloring $c_t$
of the spanning subgraph $G_t$ in the nonincreasing ear
decomposition that uses at most $\lceil\frac{n_t}{2}\rceil$ colors.
If $G$ has $2$ or $3$ ears of length 2 in the nonincreaing ear
decomposition, then $G_{t-2}$ has at most one ear of length $2$ and
$\ell(P_{t-1})=\ell(P_t)=2$. From Lemmas \ref{lem1} and \ref{lem2},
$G_{t-2}$ has a noncomplete $\lceil\frac{n_{t-2}}{2}\rceil$-rainbow
coloring $c_{t-2}$. Assume that $P_{t-1}=a_{t-1},v,b_{t-1}$ and
$P_t=a_t,v^\prime , b_t$. Since $P_{t-1}$ is a longest ear of
$G_{t-2}$, we have $a_t, b_t\in V(G_{t-2})$. Since the coloring
$c_{t-2}$ is noncomplete, $a_{t-1}$ has a noncomplete rainbow path
$P$ to one of $a_t, b_t$ (say $a_t$) such that the color $x$ in
$G_{t-2}$ does not appear in $P$. Define an edge-coloring of
$P_{t-1}$ and $P_t$ by $c(a_{t-1}v)=c(b_{t-1}v)=c(b_tv^\prime )=
x_1$ and $c(a_tv^\prime)=x$, where $x_1$ is a new color. It is clear
that $va_{t-1}Pa_tv^\prime$ is a rainbow path from $v$ to
$v^\prime$, and the obtained coloring of $G_t$ is a
$\lceil\frac{n_t}{2}\rceil$-rainbow coloring.

Now consider the case that $G$ has at least 4 ears of length 2 in
the nonincreaing ear decomposition. Suppose that
$\ell(P_{t^\prime-1})\geq 3$ and
$\ell(P_{t^\prime})=\ell(P_{t^\prime+1})=\cdots=\ell(P_t)=2$. Since
$P_i(1\leq i\leq k)$ is a longest ear, we have that $a_{t^\prime},
b_{t^\prime}, \cdots , a_t, b_t \in V(G_{t^\prime-1})$. From Lemmas
\ref{lem1} and \ref{lem2}, there exists a
$\lceil\frac{n_{t^\prime-1}}{2}\rceil$-rainbow coloring
$c_{t^\prime-1}$ of $G_{t^\prime-1}$. Color one edge of
$P_i(t^\prime \leq i\leq t)$ with $x_1$ and the other with $x_2$,
where $x_1, x_2$ are two new colors. It is obvious that $G_t$ is
rainbow connected. Since $G$ has at least $4$ ears of length $2$,
the rainbow coloring of $G_t$ uses at most
$\lceil\frac{n_t}{2}\rceil$ colors.   $\blacksquare$

From the above three lemmas and the fact that
$rc(C_n)=\lceil\frac{n}{2}\rceil \ (n\geq 4)$, we can derive our
following main result.

\begin{thm} \label{thm1} Let $G$ be a 2-connected graph of
order $n \ (n\geq 3)$. Then $rc(G)\leq\lceil\frac{n}{2}\rceil$, and
the upper bound is sharp for $n\geq 4$.
\end{thm}

Since for any two distinct vertices in a $k$-connected graph $G$ of
order $n$, there exist at least $k$ internal disjoint paths
connecting them, the diameter of $G$ is no more than
$\lfloor\frac{n}{k}\rfloor$. One could think of generalizing Theorem
\ref{thm1} to the case of higher connectivity in the obvious way,
and pose the following conjecture.

\begin{conj}\label{conj}  Let $G$ be a $k$-connected graph $G$ of order $n$.
Then $rc(G)\leq \lceil\frac{n}{k}\rceil$.
\end{conj}

\end{document}